\documentclass[reqno]{article}
\usepackage{amsmath}
\usepackage{amssymb}
\usepackage{amsthm}
\usepackage{amscd}
\usepackage{amsfonts}
\usepackage{graphicx}%
\usepackage{subfigure}
\usepackage{geometry}
\usepackage{hyperref}
\usepackage{float}
\usepackage{caption}
\usepackage[OT2,T1]{fontenc}
\usepackage[english,english]{babel}

\begin{document}

\centerline{\huge \bf On multiple translative tiling in the plane}

\bigskip\medskip
\centerline{\Large\bf Qi Yang}

\vspace{1cm}\noindent{\large\bf Abstract.}   This paper shows that in the plane, a  multiple translative tile is a multiple lattice tile.

\bigskip\noindent{\bf\large 1. Introduction}

\medskip Let $K$ be a convex body with non-empty interior int$(K)$ and boundary $\partial K$, and let $X$ be a
discrete multiset in $\mathbb{E}^{d}$. We call $K + X$ a translative tiling of $\mathbb{E}^{d}$ and call $K$ a translative tile if
$K + X =\mathbb{E}^{d}$ and any pair of translates int$(K) + \mathbf{x}_{i}$ are disjoint. In other words, $K + X$ is both a
packing and a covering in $\mathbb{E}^{d}$ . Particularly, if $X$ is a lattice in $\mathbb{E}^{d}$ , we call $K +X$ a lattice tiling of
$\mathbb{E}^{d}$ and call $K$ a lattice tile.

\smallskip Apparently, a translative tile is a convex polytope. It was shown by Minskowski \cite{8} in 1897
that every translative tile must be centrally symmetric. In 1954, Venkov \cite{11} proved that every
translative tile must be a lattice tile. Later, a new proof for this beautiful result was independently
discovered by McMullen \cite{9}.

\smallskip Let $X$ be a discrete multiset in $\mathbb{E}^{d}$ and let $k$ be a positive integer. We call $K + X$ a $k$-fold
translative tiling of $\mathbb{E}^{d}$ and call $K$ a translative $k$-tile if any point $x\in\mathbb{E}^{d}$  belongs to at least
$k$ translates of $K$ in $K + X$ and every point $x\in\mathbb{E}^{d}$ belongs to at most $k$ translates of int$(K)$
in int$(K) + X$. In other words, if $K + X$ is both a $k$-fold packing and a $k$-fold covering in $\mathbb{E}^{d}$.
Particularly, if $X$ is a lattice in $\mathbb{E}^{d}$, we call $K + X$ a $k$-fold lattice tiling of $\mathbb{E}^{d}$ and call $K$ a lattice
$k$-tile. We call $K$ a multiple translative (lattice) tile if $K$ is a translative (lattice) $k$-tile for some
positive integer $k$.

\smallskip In 1936, multiple tiling was first investigated by Furtwängler \cite{2} as a generalization of Minkowski’s conjecture on cube tiling. For more information, see \cite{10}, \cite{5} and \cite{14}.  Similar to Minkowski’s characterization, it was shown by Gravin, Robins and Shiryaev \cite{4} that a translational $k$-tile must be a centrally symmetric polytope with centrally symmetric facets.  As an analogy to the beautiful
results of Venkov \cite{11} and McMullen \cite{9}, it is natural to ask if a multiple translative tile is a multiple lattice tile.

\smallskip In 2000, Kolountzakis \cite{6} studied the structure of a multiple translative tiling by proving that,
if $D$ is a two dimensional convex domain which is not a parallelogram and $D + X$ is a multiple
tiling in $\mathbb{E}^{2}$, then $X$ must be a finite union of $2$-dimensional translated lattice. In 2013, Gravin,
Kolountzakis, Robins and Shiryae \cite{3} discovered a similar result in $\mathbb{E}^{3}$.

\smallskip Let $\tau(K)$ denote the smallest integer $k$ such that $K$ can form a $k$-fold translative tiling in $\mathbb{E}^{d}$,
and let $\tau^{*}(K)$ denote the smallest integer $k$ such that $K$ can form a $k$-fold lattice tiling in $\mathbb{E}^{d}$. For
convenience, we define $\tau(K)=\infty$ if $K$ cannot form translative tiling of any multiplicity. Clearly,
for every centrally symmetric convex polytope we have
$\tau(K)\leq\tau^{*}(K)$.

\smallskip At the end of \cite{3}, several open problems were proposed. One of them is: Prove or disprove that
if any polytope $k$ tile $\mathbb{E}^{d}$ by translations, then it is also $m$ tile $\mathbb{E}^{d}$ by lattice, for a possibly different
$m$. i.e. Prove or disprove that $\tau(K)<\infty$ imply that $\tau^{*}(K)<\infty$. This paper confirms the
two-dimensional case of this problem. We acknowledge that this result is independently discovered
by Liu \cite{7}.

\smallskip In 2017, Yang and Zong \cite{12}  studied the multiple translative tiling with given multiplicity. Later,
Zong \cite{15} characterized all the two-dimensional five-fold lattice tile. Afterwards, Yang and Zong
\cite{13}  showed that in the plane $\tau(P)=5$ imply that $\tau^{*}(P)=5$ and thus characterized all the two-dimensional translative tile. There is very little known in relations between $\tau(P)$ and $\tau^{*}(P)$, even in the plane.

\newpage\noindent The main result of  this paper is the following theorem :

\smallskip\noindent {\bf Theorem 1.} In the plane, a multiple  translative tile is a multiple lattice tile.

\vspace{1cm}\noindent{\bf\large 2. Preparation}

\medskip To prove theorem 1 we  need the following known results.

\medskip\noindent{\bf The structure of a multiple translative tiling in the plane}

\bigskip\noindent{\bf Theorem 2 (Mihail N.Kolountzakis \cite{6}).} Suppose that $P+X$ is a multiple translative tiling, where $X$ is a multiset in the plane. If $P$ is not a parallelogram, then $X$ is a finite union of two-dimensional  lattices.

\bigskip  By slightly modifying  the method used in U.Bolle \cite{1}, we can get more information  about the structure of a multiple translative tiling. Without specific statement, assume that $P+X$ is a $k$-fold translative tiling in the plane for some positive integer $k$.

\medskip\noindent{\bf Definition 1 (U.Bolle \cite{1}).} Let $L(\mathbf{e})$ be the straight line containing $\mathbf{e}$, where $\mathbf{e}$ is an edge of $P+\mathbf{x}$ and $\mathbf{x}\in X$.

 A point $\mathbf{p}\in L(\mathbf{e})$ is called a normal point if there is an $\varepsilon>0$ with
 $$(B_{\varepsilon}(\mathbf{p})\backslash{L(\mathbf{e})})\cap(\bigcup_{\mathbf{x}\in X}(\mathbf{x}+\partial{P}))=\emptyset,$$

\noindent where $B_{\varepsilon}(\mathbf{p})$ denotes the open circular disc with center $\mathbf{p}$ and radius $\varepsilon$. Since $X$ is a discrete multiset in $\mathbb{E}^{2}$, one can deduce that almost all points of $L(\mathbf{e})$ are normal and the non-normal points form a discrete set.

For normal points we define two functions $n_{i}(i=1,2)$ by

\medskip
\begin{equation}N_{i}(\mathbf{p})=\{\mathbf{x}\in X \ | \  \mathbf{ x}+P\subseteq{cl\{H_{i}\}}\  \text{and}\  \mathbf{p}\in \mathbf{x}+\partial{P}\}\end{equation}

\begin{equation}n_{i}(\mathbf{p})=|N_{i}(\mathbf{p})|\end{equation}

\noindent where $H_{i}(i=1,2)$ are the two half-planes defined by $L(\mathbf{e})$ and $ cl\{H_{i}\} $ are the closure of $\{H_{i}\}$.

In fact, $n_{1}(\mathbf{p})=n_{2}(\mathbf{p})$ for all normal points $\mathbf{p}$. For, if $B_{i}=B_{\varepsilon}(\mathbf{p})\cap H_{i}$,  each point in $B_{i}$ is cover exactly $k$ times, and if we cross $L(e)$ in $\mathbf{p}$ from $B_{1}$ to $B_{2}$, then we leave $n_{1}$ translates of $P$ and enter $n_{2}$ translates of $P$.

\medskip  Assume that $P$ is a centrally symmetric polygon with center $\mathbf{o}$ and $2m$ edges for some positive integer $m\geq 4$. Let $\mathbf{v}_{1}$,$\mathbf{v}_{2}$,...,$\mathbf{v}_{2m}$  be the 2m vertices of $P$ enumerated in the counterclockwise order. Define $\mathbf{e}_{i}=\mathbf{v}_{i+1}-\mathbf{v}_{i}$ ($1\leq i\leq 2m,\mathbf{v}_{1}=\mathbf{v}_{2m+1}$), $\mathbf{e}^{*}_{i}=\mathbf{v}_{i+m}-\mathbf{v}_{i+1}$, $1\leq i\leq m$.

\bigskip\noindent{\bf Lemma 3.} Suppose that  $\mathbf{x}\in X$. For each $i$ ($1\leq i\leq m$), either  $\mathbf{x}-\mathbf{e}_{i}$ or $\mathbf{x}-\mathbf{e}_{i}^{*}$ belongs to $X$.

\medskip\noindent{\bf Proof.}  Assume that $\mathbf{x}-\mathbf{e}_{i}^{*}\notin X$ and $\mathbf{x}-\mathbf{e}_{i}\notin X$. Let $L$ be the line  determined by $\mathbf{e}_{i}+\mathbf{x}$. Since the value of
$n_{1}$ of the normal points on the line $L$ can only change at endpoints of translates of $\mathbf{e}_{i}$.  Let $R$ be the endpoint $\mathbf{v}_{i}+\mathbf{x}$ of $\mathbf{e}_{i}+\mathbf{x}$,  let $Q_{1}$ and $Q_{2}$ be normal points of $L$ separated by $R$ such that there are only normal points between $R$ and $Q_{i}$.

Since $\mathbf{x}-\mathbf{e}_{i}\notin X$, then we have $n_{1}(Q_{1})>n_{1}(Q_{2})$ and so $n_{2}(Q_{1})\neq n_{2}(Q_{2})$ .  It is easy to see that $R$ is an endpoint of a translates of $\mathbf{e}_{i+m}$, as shown by Figure 1. By assumption, the other endpoint $R^{*}$ of $\mathbf{e}_{i}+\mathbf{x}$ cannot belong to the same translate of $\mathbf{e}_{i+m}$, then we have $n_{2}(Q_{1})< n_{2}(Q_{2})$, a contradiction. As a conclusion, we've proved the lemma 3. $\blacksquare$

\begin{figure}[H]
\centering
\includegraphics[width=2in]{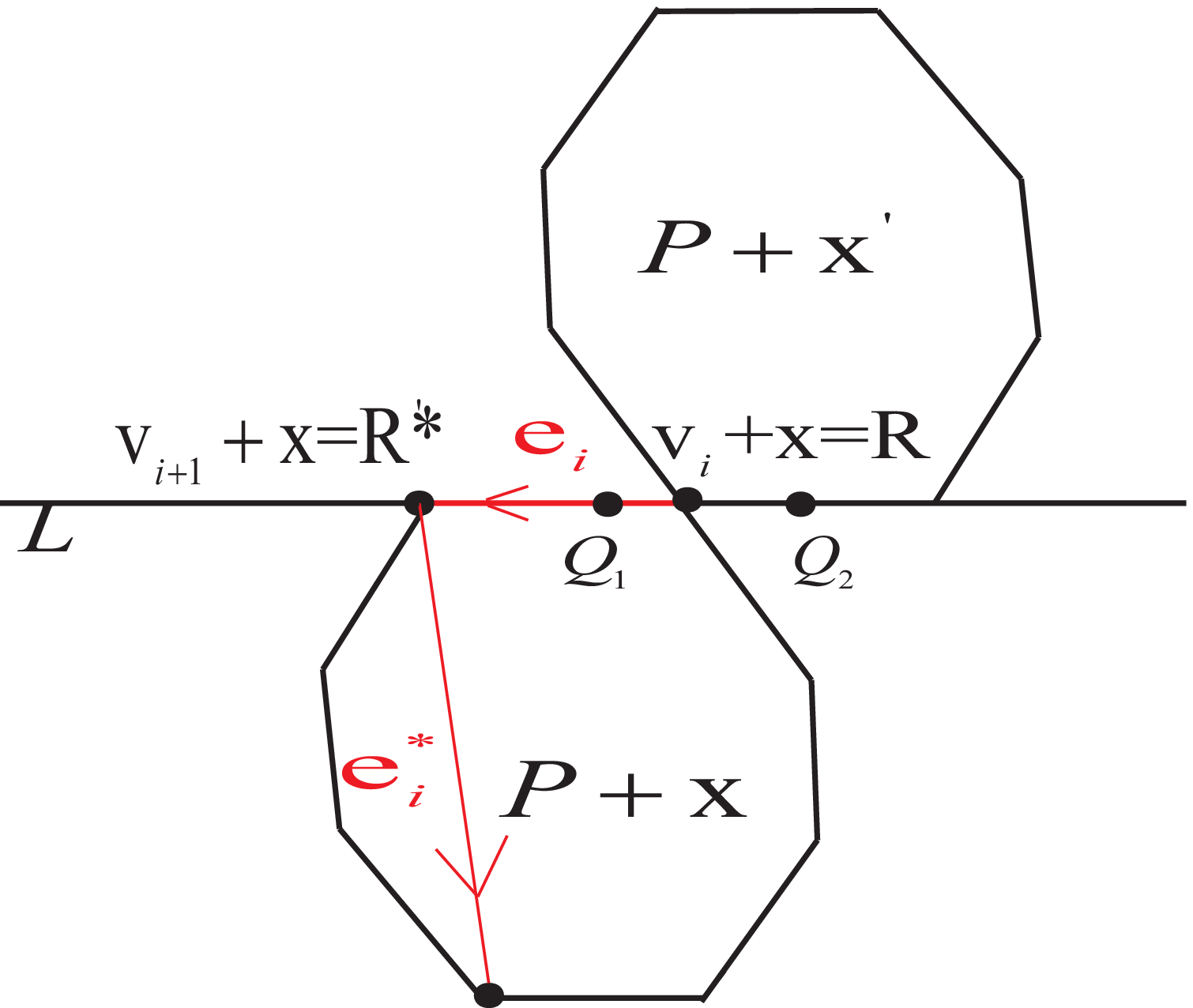}
\caption{}
\end{figure}

\noindent{\bf Theorem 4 (Bolle\cite{1}).} A convex polygon  is a $k$-fold lattice tile for a lattice $\Lambda$ and some positive integer $k$ if and only if the following conditions are satisfied:

\noindent{\bf 1.} It is centrally symmetric.

\noindent{\bf 2.} When it is centered at the origin, in the relative interior of each edge $\mathbf{g}$  there is a point of $\frac{1}{2}\Lambda$.

\noindent{\bf 3.} If the midpoint of $\mathbf{g}$ is not in $\frac{1}{2}\Lambda$ then $\mathbf{g}$ is a lattice vector of $\Lambda$.

\vspace{1cm}\noindent{\large\bf 3.  Proof of Theorem 1}

 \medskip\noindent By Theorem 2,  we can easily deduce that  $X$ can be represented as

\begin{equation}X=\bigcup\limits_{j=1}\limits^{l}\bigcup\limits_{k=1}\limits^{l_{j}}(\Lambda_{j}+\mathbf{x}^{k}_{j})(\mathbf{x}^{k}_{j}\in\mathbb{E}^{2 })\end{equation}

\noindent where $\Lambda_{j}\cap\Lambda_{j'}$ is at most one-dimensional sublattice, for any $1\leq j\neq j'\leq l$.

\bigskip\noindent{\bf Lemma 5.}  For each $1\leq j\leq l$ there exists some positive integer $\beta_{j}$ such that either $\mathbf{e}_{i}$ or  $\mathbf{e}^{*}_{i}$ belongs to $\frac{1}{\beta_{j}}\Lambda_{j}$ for all $i=1,2,...,m$.

\bigskip\noindent{\bf Proof.}   Without loss of generality, we are going to show that  there exists some positive integer $\beta_{1}$ such that either $\mathbf{e}_{i}$ or  $\mathbf{e}^{*}_{i}$ belongs to $\frac{1}{\beta_{1}}\Lambda_{1}$ for all $i=1,2,..., m$.

\medskip First, the following will be shown:

(i)  For each $1\leq i\leq m$,  \

$p_{i}\mathbf{e}_{i}+q_{i}\mathbf{e}^{*}_{i}\in\Lambda_{1}$, where  $p_{i},q_{i}(1\leq i\leq m)$ are some non-negative integers and can't be zero at the same time.

\medskip Given $\Lambda_{1}+\mathbf{x}^{k}_{1}$, the lattice $\Lambda_{1}+\mathbf{x}^{k}_{1}$ can be divided into two parts $A$ and $B$,

$$A=\{\mathbf{x}\in\Lambda_{1}+\mathbf{x}^{k}_{1}| \ \mathbf{x}-\mathbf{e}_{i}\in X \}$$

$$B=\{\mathbf{x}\in\Lambda_{1}+\mathbf{x}^{k}_{1}| \ \mathbf{x}-\mathbf{e}^{*}_{i}\in X\}$$

\noindent and obviously, $A\cup B=\Lambda_{1}+\mathbf{x}^{k}_{1}$ by lemma 3.

 Since $\Lambda_{1}\cap\Lambda_{j}$ ($j\neq1$) is at most one-dimensional sublattice, so  $((A-\mathbf{e}_{1})\cap\Lambda_{j})\cup((B-\mathbf{e}^{*}_{1})\cap\Lambda_{j})$ for $j\neq1$ is at most one-dimensional sublattice. Therefore there must be a lattice $\Lambda_{1}+\mathbf{x}^{\psi(k,1)}_{1}(1\leq\psi(k,1)\leq l_{1})$ such that

$$ (A-\mathbf{e}_{1})\cap(\Lambda_{1}+\mathbf{x}^{\psi(k,1)}_{1}) \ or \ (B-\mathbf{e}^{*}_{1})\cap(\Lambda_{1}+\mathbf{x}^{\psi(k,1)}_{1})$$

\noindent contains infinite elements and is two-dimensional.

So we  have either $\mathbf{x}^{k}_{1}-\mathbf{e}_{i}=\mathbf{x}^{\psi(k,1)}_{1}$ or $\mathbf{x}^{k}_{1}-\mathbf{e}^{*}_{i}=\mathbf{x}^{\psi(k,1)}_{1}$. Then we can define a morphism $\psi_{1}$ from $\{1,2,...,l_{1}\}$ to itself:

$$\psi_{1}(k)=\psi(k,1)$$

\noindent such that $\mathbf{x}^{k}_{1}-\mathbf{e}_{i}=\mathbf{x}^{\psi_{1}(k)}_{1}$ or $\mathbf{x}^{k}_{1}-\mathbf{e}^{*}_{i}=\mathbf{x}^{\psi_{1}(k)}_{1}$.

In this way we get   a sequence $\{a_{i}\}(a_{i}\in\{1,...,l_{1}\})$ with

 $$a_{i+1}=\psi_{1}(a_{i}) $$

Clearly, there exist two elements $a_{n_{1}},a_{n_{2}}$  of this  sequence that $a_{n_{1}}=a_{n_{2}}$, which means there exist non-negative integers  $p_{i},q_{i}$ such that

$$-p_{i} \mathbf{e}_{i}-q_{i} \mathbf{e}^{*}_{i}\in\Lambda_{1}(1\leq i\leq m)$$

\noindent where  $p_{i},q_{i}$ are not zero at same time.

Suppose that $\mathbf{u},\mathbf{v}$ are the basis vectors of $\Lambda_{1}$, then we have that

\begin{equation}p_{i}\mathbf{e}_{i}+q_{i}\mathbf{e}^{*}_{i}=a_{i}\mathbf{u}+b_{i}\mathbf{v}\end{equation}

\noindent where $a_{i},b_{i}\in\mathbb{Z}$ .

It's easy to see that

\begin{equation}\mathbf{e}^{*}_{1}=\sum^{m}_{k=2}\mathbf{e}_{k}\end{equation}

\begin{equation}\mathbf{e}^{*}_{i}=\sum^{i-1}_{k=1}-(\mathbf{e}_{k})+\sum^{m}_{k=i+1}\mathbf{e}_{k}(2\leq i\leq m)\end{equation}

\noindent So

\begin{equation}p_{1}\mathbf{e}_{1}+q_{i}\sum^{m}_{k=2}\mathbf{e}_{k}=a_{1}\mathbf{u}+b_{1}\mathbf{v}\end{equation}
\begin{equation}q_{i}\sum^{i-1}_{k=1}-(\mathbf{e}_{k})+p_{i}\mathbf{e}_{i}+q_{i}\sum^{m}_{k=i+1}\mathbf{e}_{k}=a_{i}\mathbf{u}+b_{i}\mathbf{v}(2\leq i\leq m)\end{equation}

Define $I=\{1\leq i\leq m:q_{i}=0\text{ and so}\ \mathbf{e}_{i}\in\frac{1}{p_{i}}\Lambda_{1}\}$ and $I'=\{1,2,...m\}\backslash I$, and denote the indicator function of $I'$ by $\delta_{I'}$. Then  the above equation can be represented as follows:

\begin{equation}q_{i}\sum^{i-1}_{k=1}-(\mathbf{e}_{k})\delta_{I'}(k)+\delta_{I'}(i)p_{i}\mathbf{e}_{i}+q_{i}\sum^{m}_{k=i+1}\delta_{I'}(k)\mathbf{e}_{k}=a'_{i}\mathbf{u}+b'_{i}\mathbf{v}\end{equation}

\noindent where $a'_{i},b'_{i}$ are rational numbers.

For convenience, those edges in $\{\mathbf{e}_{i}\}_{i\in I'}$ can be  re-enumerated  in original order, and denote  by $\mathbf{e}'_{j}(1\leq j\leq l')$, where $l'=|I'|$. Then we have

$$\left[\begin{matrix}
p_{1}&q_{1}&.....&q_{1}\\
-q_{2}&p_{2}&.....&q_{2}\\
...&...&...&...\\
-q_{l'}&-q_{l'}&.....&p_{l'}\\
\end{matrix}\right]
\left[\begin{matrix}
\mathbf{e}'_{1}\\
\mathbf{e}'_{2}\\
...\\
\mathbf{e}'_{l'}\\
\end{matrix}\right]
=
B
\left[\begin{matrix}
u\\
v\\
\end{matrix}\right]
$$

\noindent where all entries of $B$ are  rational numbers.

    By a series of  linear matrix transformation , we get

$$\left[\begin{matrix}
p'_{1}&1&.....&1\\
-1&p'_{2}&.....&1\\
...&...&...&...\\
-1&-1&.....&p'_{l'}\\
\end{matrix}\right]
\left[\begin{matrix}
\mathbf{e}'_{1}\\
\mathbf{e}'_{2}\\
...\\
\mathbf{e}'_{l'}\\
\end{matrix}\right]
=
B'
\left[\begin{matrix}
u\\
v\\
\end{matrix}\right]
$$

\noindent where $p'_{i}(1\leq i\leq l')$ are non-negative rational numbers and all entries of $B'$ are  rational numbers.

\medskip\noindent Define

\begin{equation}
A(p'_{1},...,p'_{l'})
=\left[\begin{matrix}
p'_{1}&1&.....&1\\
-1&p'_{2}&.....&1\\
...&...&...&...\\
-1&-1&.....&p'_{l'}\\
\end{matrix}\right]
\end{equation}

Next, we are going to prove that if   $p'_{i}\geq0$ for each $1\leq i\leq l'$ and some $p'_{j}>0$, then $A(p'_{1},...,p'_{l'})$ is invertible. We are going to prove this assertion by induction on $l'$.

It is obvious that the assertion is true while $l'=1$. Assume that the assertion holds for  $l'\leq n-1$.

\medskip Next, we are  going to prove this assertion holds when $l'=n$.  The determinant of $A(p'_{1},...,p'_{l'})$  can be seen as a function of variables $p'_{1},...,p'_{l'}$.  It is easy to calculate that the partial derivative of $A(p'_{1},...,p'_{l'})$ respect to $p'_{i}(1\leq i\leq l')$ is $A(p'_{1},...,p'_{i-1},p'_{i+1},...,p'_{l'})$, by the assumption, which is greater than 0.

 Clearly, $A(0,0,...0)=0$ when $l'$ is odd, otherwise $A(0,0,...0)=1$. As a conclusion, we've proved the assertion.

  When $p'_{i}=0$ for all $1\leq i\leq l'$, then $q_{i}\mathbf{e}^{*}_{i}\in\Lambda_{1}$ for $i\in I'$ and  $\mathbf{e}_{i}\in\Lambda_{1}$ for each $i\in I$.

  When $p'_{i}\neq0$ for some $1\leq i\leq l'$, then $A(p'_{1},...,p'_{l'})$ is invertible, so there exists an invertible matrix $A(p'_{1},...,p'_{l'})^{-1}$ such that

$$
\left[\begin{matrix}
\mathbf{e}'_{1}\\
\mathbf{e}'_{2}\\
...\\
\mathbf{e}'_{l}\\
\end{matrix}\right]
=
A(p'_{1},...,p'_{l'})^{-1}
B'
\left[\begin{matrix}
\mathbf{u}\\
\mathbf{v}\\
\end{matrix}\right]
$$

  Since  all  entries of $A(p'_{1},...,p'_{l'})^{-1}$  and $B'$ are  rational numbers, there exists a positive integer $\beta_{1}$ such that  all entries of $\beta_{1}A(p'_{1},...,p'_{l'})^{-1}B'$ are integers, i.e.  $\mathbf{e}_{i}\in\frac{1}{\beta_{1}}\Lambda_{1}$ for each $i\in I'$.

 \medskip As a conclusion, we've proved that there exist $\beta_{1}$ such that for all $1\leq i\leq m$, either $\mathbf{e}_{i}$ or  $\mathbf{e}^{*}_{i}$ belongs to $\frac{1}{\beta_{1}}\Lambda_{1}$. Similarly, we can prove the case of any other lattice $\Lambda_{j}(1<j\leq l)$, lemma 5 is proved. $\blacksquare$

\vspace{0.7cm} By lemma  4 and lemma 5, in order to prove Theorem 1, it is suffice to show that for some $\Lambda_{j}$, there is a point of $\frac{1}{2\beta_{j}}\Lambda_{j}$ in the relative interior of $\mathbf{e}_{i}$ for each $i(1\leq i\leq m)$.

\medskip\noindent{\bf Lemma 6.} ii) With above notations, for some $j(1\leq j\leq l)$, there  exists a point of $\frac{1}{2\beta_{j}}\Lambda_{j}$ in the relative interior of $\mathbf{e}_{i}$ for each $1\leq i\leq m$, where $\beta_{j}$  is some positive integer.

\medskip\noindent{\bf Proof.}
To prove lemma 6, it is sufficient to discuss the following two cases:

\medskip \noindent {\bf Case 1.} For  some $\Lambda_{j}$, there exists some $i(1\leq i\leq m)$ such that $p\mathbf{e}_{i}$ and $q \mathbf{e}^{*}_{i}\in\Lambda_{j}$, where $p, q\in\mathbb{Q}\backslash{\{0\}}$.

\medskip By lemma 5, it is easy to deduce that for all  $i(1\leq i\leq m)$, $q_{i}\mathbf{e}^{*}_{i}\in\Lambda_{i}$ for some positive $q_{i}\in\mathbb{Q}$, so  there exists some positive integer $\beta_{j}$ such that $\mathbf{e}^{*}_{i}\in\frac{1}{\beta_{j}}\Lambda_{j}$ for all $\mathbf{e}^{*}_{i}$, the middle point of every edge $\mathbf{e}_{i}$ of $P$ belongs to $\frac{1}{2\beta_{j}}\Lambda_{j}$, lemma 6 holds.

\vspace{0.7cm}\noindent {\bf Case 2.} For each $\Lambda_{j}(1\leq j\leq l)$, $p \mathbf{e}_{i}$ and $q \mathbf{e}^{*}_{i}$  can not belong to $\Lambda_{j}$ at the same time for each $1\leq i\leq m$, where $p, q$ are arbitrary non-zero rational numbers.

\medskip In this case, we are going to prove that lemma 6 holds for $\Lambda_{1}$ ( The argument for other lattice $\Lambda_{j}$ is similar).  If $\mathbf{e}^{*}_{i}$ belongs to $\frac{1}{\beta} \Lambda_{1}$ for some positive integer $\beta$, then the midpoint of $\mathbf{e}_{i}$ belongs to the lattice $\frac{1}{2\beta} \Lambda_{1}$.  Otherwise $\mathbf{e}_{i}\in\frac{1}{\beta} \Lambda_{1}$ for some positive integer $\beta$, then it is suffice  to show that there exist  $\mathbf{g}\in\frac{1}{\beta} \Lambda_{1}$ such that $\mathbf{g}-\mathbf{e}^{*}_{i}=\lambda \mathbf{e}_{i}$ for some $\lambda\in\mathbb{R}$.


 \bigskip Assume that $\mathbf{e}^{*}_{1}\notin\frac{1}{\beta}\Lambda_{1}$ for any positive integer $\beta$. Then we have $\mathbf{e}_{1}\in\frac{1}{\beta}\Lambda_{1}$ for some positive integer $\beta$, which means that $\mu\mathbf{e}_{1}\in\Lambda_{1}$ for some positive rational number $\mu$. Define  $\Omega=\{\Lambda_{j}+\mathbf{x}^{k}_{j},1\leq j\leq l,1\leq k\leq l_{j}\}$.  For convenience, we may assume that $\mu\mathbf{e}_{1}$ and $\mathbf{u}$ are the basis vectors of lattice $\Lambda_{1}$. Denote the line containing $\mathbf{e}_{1}+\mathbf{x^{1}_{1}}+n\mathbf{u}(n\in\mathbb{Z})$ by $L(\mathbf{x}^{1}_{1},n)$.

\medskip   Let $\mathbf{p}^{1,n}_{j}$ be  a normal point in the relative interior of $\mathbf{e}_{1}+\mathbf{x^{1}_{1}}+j\mu\mathbf{e}_{1}+n\mathbf{u}$ for each $j\in\mathbb{Z}$. Since  $n_{1}(\mathbf{p}^{1,n}_{j})=n_{2}(\mathbf{p}^{1,n}_{j})>0$, then we  define

$$\Gamma_{1,n}=\{\mathbf{y}^{1,n}_{j}:\mathbf{y}^{1,n}_{j}\in N_{2}(\mathbf{p}^{1,n}_{j})\}_{j\in\mathbb{Z}}$$

 Apparently, we can find an infinite subset $\Gamma'_{1,n}$ of  $\Gamma_{1,n}$ which its elements are contained  in the same translated lattice in $\Omega$, denote the corresponding lattice by $\Lambda(\Gamma_{1,n})$.

\medskip Define $W_{1}=\{n\in\mathbb{N}:p\mathbf{e}_{1}\notin\Lambda(\Gamma_{1,n}) \text{for any}\ p\in\mathbb{Q}\}$. If $|W_{1}|<\infty$, then $|\mathbb{N}\backslash W_{1}|=\infty$ and so we can find an infinite subset $W'_{1}$ of $\mathbb{N}\backslash W_{1}$ that   the corresponding translated lattice $\Lambda(\Gamma_{1,n})$ for every $n\in W'_{1}$ is identical, denote this lattice by $\Lambda(1)$. For each $n\in W'_{1}$, let $\mathbf{p}^{2,n}_{j}$ be a normal point in the relative interior of the edge $\mathbf{e}_{1}+\mathbf{y}^{1,n}_{j}$ for $\mathbf{y}^{1,n}_{j}\in\Gamma'_{1,n}$. Since $n_{1}(\mathbf{p}^{2,n}_{j})=n_{2}(\mathbf{p}^{2,n}_{j})>0$, by analogy, define $\Gamma_{2,n}$ as follows:

\begin{equation}\Gamma_{2,n}=\{\mathbf{y}^{2,n}_{j}:\mathbf{y}^{2,n}_{j}\in N_{2}(\mathbf{p}^{2,n}_{j})\ \text{for}\ j \ \text{such that }\mathbf{y}^{1,n}_{j}\in\Gamma'_{1,n}\}\end{equation}

\noindent then we can  find an infinite subset $\Gamma'_{2,n}$ of  $\Gamma_{2,n}$ such that its elements are contained  in the same translated lattice in $\Omega$, denote the corresponding lattice  by $\Lambda(\Gamma_{2,n})$.

  \begin{equation}\Gamma'_{2,n}=\{\mathbf{y}^{2,n}_{j}:\mathbf{y}_{j,n}\in\Lambda(\Gamma_{2,n})\}\end{equation}

\medskip
\begin{figure}[H]

\centering

\includegraphics[width=2.5in]{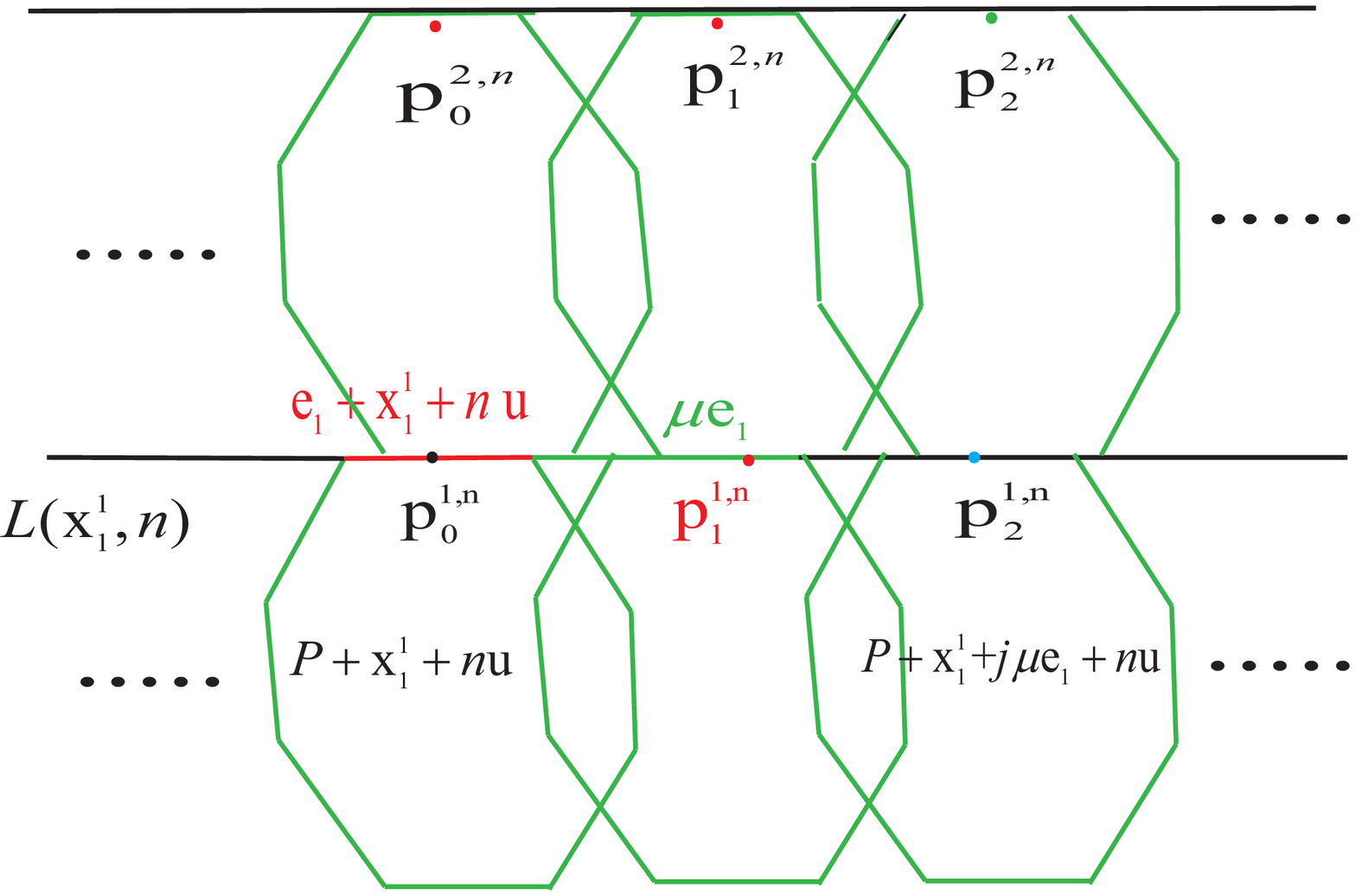}

\caption{}
\end{figure}

   Similarly, we can define the set  $W_{2}=\{n\in W'_{1}:p\mathbf{e}_{1}\notin\Lambda(\Gamma_{2,n}) \text{for any}\ p\in\mathbb{Q}\}$. If $|W_{2}|<\infty$, then we can define an infinite subset $W'_{2}$ of $W'_{1}\backslash W_{2}$. Next, we are going to use induction to define the notations appeared in the following argument.

   \medskip

   Suppose that we have $|W_{k}|<\infty(k\geq 2)$, and $|W'_{k-1}|=\infty$. Then there is an infinite subset $W'_{k}$ of $W'_{k-1}\backslash W_{k}$ that   the corresponding translated lattice $\Lambda(\Gamma_{k,n})$ for every $n\in W'_{k}$ is identical,  denote this lattice by $\Lambda(k)$.

  \medskip

   Let $\mathbf{p}^{k+1,n}_{j}$ be a normal point in the relative interior of the edge $\mathbf{e}_{1}+\mathbf{y}^{k,n}_{j}$ for $\mathbf{y}^{k,n}_{j}\in\Gamma'_{k,n}$ and $n\in W'_{k}$.

  \medskip

  Since $n_{1}(\mathbf{p}^{k+1,n}_{j})=n_{2}(\mathbf{p}^{k+1,n}_{j})>0$,  define $\Gamma_{k+1,n}$ as follows:

\begin{equation}\Gamma_{k+1,n}=\{\mathbf{y}^{k+1,n}_{j}:\mathbf{y}^{k+1,n}_{j}\in N_{2}(\mathbf{p}^{k+1,n}_{j})\ \text{for}\ j \ \text{such that }\mathbf{y}^{k,n}_{j}\in\Gamma'_{k,n}\}\end{equation}

\noindent then we can  find an infinite subset $\Gamma'_{k+1,n}$ of  $\Gamma_{k+1,n}$ which its elements are contained  in the same translated lattice in $\Omega$,  denote the corresponding lattice denote by $\Lambda(\Gamma_{k+1,n})$.

  $$\Gamma'_{k+1,n}=\{\mathbf{y}^{k+1,n}_{j}:\mathbf{y}^{k+1,n}_{j}\in\Lambda(\Gamma_{k+1,n})\cap\Gamma'_{k+1,n} \}$$

\noindent  and define $W_{k+1}=\{n\in W'_{k}:p\mathbf{e}_{1}\notin\Lambda(k+1,n)\ \text{for any}\ p\in\mathbb{Q}\}$.

\bigskip Since $|\Omega|<\infty$, there are only two possible cases: \textbf{(a)} $|W_{k^{*}}|=\infty$ for some positive integer $k^{*}(k^{*}\geq1)$ ; \textbf{ (b)} $\Lambda(k')=\Lambda(k'')$ for two different positive integers $k',k''$.

\medskip\noindent{\large\bf (a)} If $|W_{k^{*}}|=\infty$ for some positive integers $k^{*}(k^{*}\geq 1)$, then there exist two elements $n_{1},n_{2}$ of $W_{k^{*}}$ such that $\Lambda(k^{*},n_{1})=\Lambda(k^{*},n_{2})$.

By the definition of $\Lambda(k^{*},n_{1})$, $\mu'*\mathbf{e}_{1}\in\Lambda(k^{*},n_{1})$ for some $\mu'\in\mathbb{R}$, and by the definition of  $W_{k^{*}}$, $\mu'$ is an irrational number. Again, by the definition of $W_{k^{*}}$, $q*\mathbf{e}^{*}_{1}\in\Lambda(k^{*},n_{1})$ for some positive rational number $q$.

 By the definition of $\Lambda(k^{*},n_{1})$, there exist  two elements $\mathbf{y}^{k^{*},n_{1}}_{j'}\in\Lambda(k^{*},n_{1})$ and $\mathbf{y}^{k^{*},n_{2}}_{j''}\in\Lambda(k^{*},n_{1})$ such that:

\begin{equation}\mathbf{y}^{k^{*},n_{1}}_{j'}=\mathbf{x}^{1}_{1}+n_{1}\mathbf{u}-k^{*}\mathbf{e}^{*}_{1}+\lambda_{3}\mathbf{e}_{1}\end{equation}

\begin{equation}\mathbf{y}^{k^{*},n_{2}}_{j''}=\mathbf{x}^{1}_{1}++n_{2}\mathbf{u}-k^{*}\mathbf{e}^{*}_{1}+\lambda_{4}\mathbf{e}_{1}\end{equation}

\noindent where $\lambda_{3},\lambda_{4}\in\mathbb{R}$ and some  $j',j''\in\mathbb{Z}$.

Combined (14)-(15),

\begin{equation}\mathbf{y}^{k^{*},n_{1}}_{j'}-\mathbf{y}^{k^{*},n_{2}}_{j''}=(n_{1}-n_{2})\mathbf{u}+(\lambda_{3}-\lambda_{4})\mathbf{e}_{1}\end{equation}

Since $\mu'\mathbf{e}_{1}$ and $q\mathbf{e}^{*}_{1}$ are linearly independent, then we have

\begin{equation}\mathbf{y}^{k^{*},n_{1}}_{j'}-\mathbf{y}^{k^{*},n_{2}}_{j''}=a\mu'\mathbf{e}_{1}+bq\mathbf{e}^{*}_{1}\end{equation}

\noindent where $a,b\in\mathbb{Q}$ and $bq\neq0$ (since $n_{1}\neq n_{2}$).

 Since $\mathbf{u}$ and $\mathbf{e}_{1}$ are linearly independent, then we have

\begin{equation}\mathbf{e}^{*}_{1}=\lambda_{5}\mathbf{u}+\lambda_{6}\mathbf{e}_{1}\end{equation}

   When $\lambda_{5}\in\mathbb{Q}$, then there is a point $\mathbf{g}$ in the lattice $\frac{1}{\beta}\Lambda_{1}$ such that

$$\mathbf{g}-\mathbf{e}^{*}_{1}=\lambda_{6}\mathbf{e}_{1}$$

\noindent which means that for some positive integer $\beta$, there is a lattice point of  $\frac{1}{\beta}\Lambda_{1}$ in the relative interior of $\mathbf{e}_{1}$.

\medskip When $\lambda_{5}\notin\mathbb{Q}$, combined (16),(17) and (18), we have

\begin{equation}\mathbf{y}^{k^{*},n_{1}}_{j'}-\mathbf{y}^{k^{*},n_{2}}_{j''}=a\mu'\mathbf{e}_{1}+bq(\lambda_{5}\mathbf{u}+\lambda_{6}\mathbf{e}_{1})\end{equation}

\noindent and

\begin{equation}(n_{1}-n_{2})\mathbf{u}+(\lambda_{3}-\lambda_{4})\mathbf{e}_{1}=a\mu'\mathbf{e}_{1}+bq(\lambda_{5}\mathbf{u}+\lambda_{6}\mathbf{e}_{1})\end{equation}

\begin{equation}(n_{1}-n_{2}-bq\lambda_{5})\mathbf{u}=(a\mu'+bq\lambda_{6}+\lambda_{4}-\lambda_{3})\mathbf{e}_{1}\end{equation}

 Since $\lambda_{5}\notin\mathbb{Q}$ , the left side of the equation (22) is not equal to zero. But $\mathbf{u}$ and $\mathbf{e}_{1}$ is linearly independent over $\mathbb{R}$, then we get a contradiction.

\medskip\noindent {\bf (b)} $\Lambda(k')=\Lambda(k'')$. Suppose that $k'<k''$. By the definition of $W'(k')$ and $W'(k'')$, we have that $W'(k'')\subset W'(k')$.  Let $n_{1}\in W'(k'')$ , let $\mathbf{y}^{k',n_{1}}_{j'}\in\Lambda(k'')$  and $\mathbf{y}^{k'',n_{1}}_{j''}\in\Lambda(k'')$ for some positive integers $j',j''$,

\begin{equation}\mathbf{y}^{k',n_{1}}_{j'}=\mathbf{x}^{1}_{1}+n_{1}\mathbf{u}-k'\mathbf{e}^{*}_{1}+\lambda_{3}\mathbf{e}_{1}\end{equation}

\begin{equation}\mathbf{y}^{k'',n_{1}}_{j''}=\mathbf{x}^{1}_{1}+n_{1}\mathbf{u}-k''\mathbf{e}^{*}_{1}+\lambda_{4}\mathbf{e}_{1}\end{equation}

\noindent for some real numbers $\lambda_{3},\lambda_{4}$.

Then we have

\begin{equation}\mathbf{y}^{k'',n_{1}}_{j''}-\mathbf{y}^{k',n_{1}}_{j'}=(k''-k')\mathbf{e}^{*}_{1}+(\lambda_{3}-\lambda_{4})\mathbf{e}_{1}.\end{equation}

\medskip Let $n_{2}\in W'(k')$ and $n_{1}\neq n_{2}$, and let $\mathbf{y}^{k',n_{2}}_{j*}\in\Lambda(k')$,  then we have

\begin{equation}\mathbf{y}^{k',n_{2}}_{j*}=\mathbf{x}^{1}_{1}+n_{2}\mathbf{u}-k'\mathbf{e}^{*}_{1}+\lambda_{5}\mathbf{e}_{1}\end{equation}

\begin{equation}\mathbf{y}^{k',n_{2}}_{j*}-\mathbf{y}^{k',n_{1}}_{j'}=(n_{2}-n_{1})\mathbf{u}+(\lambda_{5}-\lambda_{3})\mathbf{e}_{1}\end{equation}

\medskip Suppose that $\frac{1}{h}(\mathbf{y}^{k',n_{2}}_{j*}-\mathbf{y}^{k',n_{1}}_{j'})\in\Lambda(k')$ for some positive integer $h$, by the definition of $\Lambda(k')$, $\mu\mathbf{e}_{1}\in\Lambda(k')$ for some real number $\mu$. Since  $\Lambda(k')=\Lambda(k'')$,     $\mathbf{y}^{k'',n_{1}}_{j''}-\mathbf{y}^{k',n_{1}}_{j'}$ can be represented as the linear combination of $\frac{1}{h}(\mathbf{y}^{k',n_{2}}_{j*}-\mathbf{y}^{k',n_{1}}_{j'})$ and $\mu\mathbf{e}_{1}$:

\begin{equation}\mathbf{y}^{k'',n_{1}}_{j''}-\mathbf{y}^{k',n_{1}}_{j'}=\frac{z_{1}}{h}(\mathbf{y}^{k',n_{2}}_{j*}-\mathbf{y}^{k',n_{1}}_{j'})+z_{2}\mu\mathbf{e}_{1}\end{equation}

\noindent where $z_{1},z_{2}\in\mathbb{Q}$.

\medskip Combined equations (24),(26) and (27), we have

\begin{equation}(k''-k')\mathbf{e}^{*}_{1}+(\lambda_{3}-\lambda_{4})\mathbf{e}_{1}=z_{1}*((n_{2}-n_{1})\mathbf{u}+\frac{1}{h}(\lambda_{5}-\lambda_{3})\mathbf{e}_{1})+z_{2}\mu\mathbf{e}_{1}\end{equation}

\noindent Simplify the equation (29), we have

$$\mathbf{e}^{*}_{1}=q_{1}\mathbf{u}+\mu'\mathbf{e}_{1}$$

\noindent for some  rational number $q_{1}$ and some real number $\mu'$.

 As a conclusion, we've proved that there is a positive integer $\beta$ such that there is a lattice point $\mathbf{g}$ in the lattice $\frac{1}{\beta}\Lambda_{1}$ such that

$$\mathbf{g}-\mathbf{e}^{*}_{1}=\lambda\mathbf{e}_{1}$$

\noindent for some real number $\lambda$.

\medskip Similarly, for each $i (1\leq i\leq m)$, we can prove that there is a positive integer $\beta_{1}$ such that there is a lattice point $g$ in the lattice $\frac{1}{\beta_{1}}\Lambda_{1}$ that

$$g-\mathbf{e}^{*}_{i}=\lambda_{i}\mathbf{e}_{i}$$

\noindent for some real number $\lambda_{i}$. As a conclusion of the above cases, we've proved lemma 6.

 By Theorem 4, lemma 5  and 6, we've prove that there exist some $j$, for some positive integer $\beta_{j}$, $P+\frac{1}{\beta_{j}}\Lambda_{j}$ is a multiple lattice tiling, Theorem 1 is proved.

\vspace{0.5cm}\noindent{\bf Acknowledgements.}  For helpful comments and suggestions, the author is grateful to Professor C.Zong. This work is supported by 973 Program 2013CB834201.

\bibliographystyle{amsplain}
{}

\vspace{0.5cm}\noindent Qi Yang, School of Mathematical Science, Peking University, Beijing 100871. China

\bf\noindent\textit{Email: Yangqi07@pku.edu.cn}
\end{document}